\documentclass[9pt,amstex]{article}
\usepackage{amssymb}
\usepackage{latexsym}

\newtheorem{dfn}{Definition}[section] 
\newtheorem{rmk}{Remark}[section]
\newtheorem{thm}{Theorem}[section] 

\newtheorem{prop}{Proposition}[section] 
\newtheorem{lem}{Lemma}[section]

\newtheorem{ex}{Example}[section] 
\newcommand{\Pf}{{\em Proof}. }
\newcommand{\EPf}
{%
\mbox{}%
\nolinebreak%
\hfill%
\rule{2mm}{2mm}%
\medbreak%
\par%
}
\newcommand{\C}{\mathbb C} 

\newcommand{\R}{\mathbb R}

\newcommand{\EE}{{\cal E}_{\hbar}}
\newcommand{\e}{{\cal E}}
\newcommand{\sech}{\mbox{sech}} 
\newcommand{\arcsinh}{\mbox{\rm arcsinh}} 
\newcommand{\g}{{\cal G}{}} 
\renewcommand{\u}{{\cal U}{}} 

\renewcommand{\k}{{\cal K}{}} 

\newcommand{\p}{{\cal P}{}} 
 
\newcommand{\iw}{{\cal I}{}} 
\newcommand{\n}{{\cal N}{}} 
 
\newcommand{\s}{{\cal S}{}} 
\newcommand{\ddto}{{\frac{d}{dt}|_{0}}{}}
\renewcommand{\j}{{\cal J}{}} 
\renewcommand{\o}{{\cal O}{}} 
\renewcommand{\r}{{\cal R}{}}

\renewcommand{\a}{{\cal A}{}}

\newcommand{\z}{{\cal Z}{}}

\renewcommand{\k}{{\cal K}{}}

\begin{document}
\title{Strict Deformation Quantization for Actions of a Class of Symplectic Lie Groups}
\author{{\bf Pierre Bieliavsky}\\
Universit\'e Libre de Bruxelles\\
Belgium\\
and\\
{\bf Marc Massar}\\
Vrije Universiteit van Brussel\\
Belgium}
\maketitle
\begin{abstract}
We present explicit universal strict deformation quantization formulae for 
actions of Iwasawa subgroups $AN$ of $SU(1,n)$. This answers a 
question raised by Rieffel in \cite{R2}.
\end{abstract}
\section*{Introduction}
In \cite{R1}, Rieffel describes explicitly strict 
quantization for actions of $\R^{d}$. Roughly, the idea is as follows.
Let $A$ be an associative (topological) algebra and let $\alpha:\R^d\times A\to A$ be an 
action of $\R^d$ on $A$ by automorphisms.
Then, if the situation is regular enough, one can give a sense to the 
following product~:
\begin{equation}\label{RIEFFEL}
a\star_Jb\stackrel{\mbox{def.}}{=}\int\alpha_x(a)\alpha_{Jy}(b)\, e^{i<x,y>}dx\, dy
\end{equation}
as an oscillatory integral. In the previous formula, $a$ and $b$ are elements 
of $A$, $<\, , \,>$ is a Euclidean scalar product on $\R^d$ and $J$ is 
a skewsymmetric matrix $J\in so(\R^d,\,<,>)$. Also, $dx\, dy$ denotes
a Haar measure on $\R^d\times\R^d$. The product $\star_J$ 
is a deformation of the product on $A$. Indeed, when $J$ tends to $0$ in 
$so(\R^d,\,<,>)$, $\star_J$ tends to the product on $A$. Moreover, the 
product $\star_J$ turns out to be associative. In the case $M$ is some 
manifold $\R^d $ acts on by diffeomorphisms, one gets an action on 
$A=C_\infty(M)$. Therefore, introducing a deformation parameter $\hbar$ multiplying $J$,
 the above construction yields a deformation 
quantization of the Poisson bracket on $M$ induced by the action and the 
data of $J$.
Formula (\ref{RIEFFEL}) is universal in the sense that it is valid for any 
action of $\R^d$. \\
\noindent An approach to universal deformation formulae (UDF) has 
been proposed by Giaquinto and Zhang \cite{GZ} within the formal 
framework i.e. the resulting deformation is a formal power series, 
such as a star product for instance. Among other things, Giaquinto 
and Zhang present there a beautiful explicit formal UDF for actions
of the (non-Abelian) Lie group ``$ax+b$'' of affine transformations 
of the real line. As observed by Rieffel in \cite{R2}, this suggests 
that, at least for the group ``$ax+b$'', there should exist UDF's at 
the analytical (strict) level.

It has been observed in \cite{Bi2} that a possible approach to ``universal" 
deformations for (non Abelian) group actions is to study a particular class of 
three-point kernels on group manifolds. Roughly speaking, this means 
the following. 
Let $G$ be a Lie group endowed with a left-invariant Haar 
measure $\mu$. Assume the existence of a (non trivial) function space $\e\subset C(G)$ and 
a three-point kernel $K\in C(G\times G\times G)$ such that
\begin{enumerate}
\item[(i)] $K$ is invariant under the diagonal left action of $G$ on $G\times G\times G$;
\item[(ii)] for all $u,v\in\e$, the formula 
$$
u\star^Gv(x)=\int_{G\times G}u(g)\,v(h)\, K(x,g,h)\, \mu_g\mu_h
$$
defines an element $u\star^Gv$ of $\e$ \footnote{This condition 
is too strong in general. The product formula should hold only on a 
dense subset of $\e$, provided a topological framework is defined 
(see e.g. \cite{R1} or Theorem \ref{THM}).};
\item[(iii)] $(\e,\star^G)$ is an associative algebra.
\end{enumerate}
Then, if an action $\alpha:G\times A\to A$ of $G$ on some (topological) algebra $A$
is given, one checks that formally the formula
\begin{equation}\label{UNIV}
a\star^Ab\stackrel{\mbox{def.}}{=}\int_{G\times 
G}K(e,g,h)\alpha_g(a)\,\alpha_h(b)\, \mu_g\mu_h,
\end{equation}
where $e$ is the unit element of $G$, defines an associative product on $A$, provided 
some regularity condition is fulfilled (e.g. the function $G\to A:g\to\alpha_g(a)$ has the same type of 
regularity as the elements of $\e$).\\
For example, in the case of the Abelian group $G=\R^d$, Rieffel's product (\ref{RIEFFEL}) 
coincides with (\ref{UNIV}) when
$$
K(x,y,z)=\exp\left(\frac{2i}{\hbar}\left\{<x,Jy>+<y,Jz>+<z,Jx>\right\}\right)
$$
and when $\e=C_\infty(\R^d)$.

In this paper, we explicitly describe universal deformation three-point kernels $K$ as above
for a class of solvable Lie groups (see Theorem \ref{THM}). Those are Iwasawa $AN$ subgroups 
of $SU(1,n)$. In particular, any action of $SU(1,n)$ yields an action 
of $AN$. This should 
provide an interesting class of strict deformation quantizations of (non regular) 
Poisson manifolds such as quantum flag manifolds. However, this last point will not be 
investigated in the present article (see nevertheless Section \ref{Schubert}).\\
This note is organized as follows.
\begin{enumerate}
\item[\ref{Buxtehude}]{\bf Rank one non-compact Hermitian symmetric spaces}\\
Let $G=SU(1,n)$ and $R=AN$ be an Iwasawa subgroup, i.e. $R$ is the 
$AN$ factor in an Iwasawa decomposition $G=ANK$ of $G$.  Then the group manifold $R$ is 
equivariantly diffeomorphic to the Hermitian symmetric space $M=G/K$. In 
particular, the Lie group $R$ carries an $R$-left-invariant symplectic structure $\omega$.
In this section, we describe a global Darboux chart on $(R,\omega)$ which 
will be important later on.
\item[\ref{Scarlatti}]{\bf Guessing the deformed product formula}\\
The proof of our main result (Theorem \ref{THM}) does not depend on the present section.
Nevertheless, it explains {\sl how} Formula (\ref{PRODUCT}) has been found.
\item[\ref{Vivaldi}]{\bf $AN$-covariant Moyal star products}\\
The main property of our Darboux chart on $(R,\omega)$ (cf. Section 
\ref{Buxtehude}) 
is that, with respect to this chart, the formal Moyal star product is covariant in Arnal's sense 
(see Definition \ref{HAM}). One can therefore apply techniques of 
star-representation theory 
\cite{BFFLS, AC, Fr}
to analyze the star-representation of $R$ arising from the covariance property.
\item[\ref{Haendel}]{\bf The $\z_\hbar$ transform}\\
Using the cocycle defining the star representation, we introduce some 
kind of Fourier integral operator (the $\z_\hbar$ transform) which intertwines the 
pointwise commutative product on $C^\infty(R)[[\hbar]]$ with an $\hbar$-dependent (commutative)
product with respect to which the Lie algebra of $R$ acts by derivations. The 
``commutative manifold'' underlying the latter product therefore carries an $R$-invariant 
deformation quantization. This $R$-space turns out to be formally $R$-equivariantly
isomorphic to the group $R$.
\item[\ref{Bach}]{\bf Strict Quantization}\\
We define a one parameter family of (Fr\'echet) function spaces $\{\e_\hbar\}_{\hbar\geq 0}$, 
each of them endowed with an associative product $\star_\hbar$, such that
\begin{enumerate}
\item[(i)] for all $\hbar\geq 0$, one has a dense inclusion
$$
C^\infty_c(R)\subset\e_\hbar;
$$
\item[(ii)] for all $u,v\in C_{c}^\infty(R)$, the product reads
$$
u\star_\hbar v(x)=\hbar^{-2\dim R}\int_{R\times R} e^{\frac{i}{\hbar}S(x,y,z)}a(x,y,z)\, u(y)\, v(z)\,
dy\, dz
$$
where $K=a\,e^{\frac{i}{\hbar}S}$ is an $R$-left-invariant three-point kernel on $R\times R
\times R$, and where $dy\, dz$ is a left invariant Haar measure on $R\times R$ (see Theorem \ref{THM}).
\end{enumerate}
\item[\ref{Mozart}]{\bf The case $n=1$}\\
We relate the present construction in the case $n=1$ with a strict quantization of the 
symplectic symmetric space $(SO(1,1)\times\R^2)/\R$ obtained in 
\cite{Bi2}. 
This leads to a relation between our work and 
Unterberger's 
composition formula for the one-dimensional Klein-Gordon calculus 
\cite{U1, U2}.
\item[\ref{Schubert}]{\bf Remark for further developments}\\
\end{enumerate}

\noindent{\bf Acknowledgments}\\
We thank Marc Rieffel for usefull comments.
The first author has been partially supported by the Communaut\'e Fran\c caise de 
Belgique, through an Action de Recherche Concert\'ee de la Direction 
de la Recherche Scientifique. The second author is supported by the 
European Commission RTN programme HPRN-CT-2000-00131 in which he is 
associated with K.U. Leuven.
\section{Rank one non-compact Hermitian symmetric 
spaces}\label{Buxtehude}
\subsection{Notations and preliminaries} 
We refer to \cite{H} for 
the general theory of (Riemannian) symmetric spaces.\\
Let $G$ be a connected real {\bf simple Lie group} with Lie algebra $\g$. 
We denote by $B$ the {\bf Killing form} on $\g$. Let 
$\sigma$ be a Cartan involution of $\g$ with Cartan decomposition 
$$
\g=\k\oplus\p \quad (\,\sigma=id_{\k}\oplus(-id_{\p})\,).
$$
Fix a maximal Abelian subalgebra $\a$ contained in $\p$. The dimension of 
$\a$  is called the (real) {\bf rank} of $\g$. Let $\Phi$ be the {\bf root 
system} of $\g$ with respect to $\a$ and fix a positive root system $\Phi^+ 
\quad(\Phi=\Phi^+\cup(-\Phi^+))$. Denote by $\g_\alpha$ the {\bf weight 
space} corresponding to $\alpha\in\Phi\cup\{0\}$. One then has
\begin{equation}\label{BOREL}
\g=\sigma(\n)\oplus\g_0\oplus\n
\end{equation}
with
$$
\n=\sum_{\alpha\in\Phi^+}\g_{\alpha}\mbox{ and } 
\g_{0}=\z_{\k}(\a)\oplus\a
$$
where $\z_{\k}(\a)$ denotes the centralizer of $\a$ in $\k$. One also 
has the {\bf Iwasawa decomposition}
$$
\g=\k\oplus\a\oplus\n
$$
which induces a global analytic diffeomorphism~:
$$
A\times N\times K \to G~: (a,n,k)\to ank
$$
where $K$ (resp. $A$, resp. $N$) denotes the (connected) analytic 
subgroup of $G$ with algebra $\k$ (resp. $\a$, resp. $\n$). The 
Iwasawa group decomposition therefore induces a global 
diffeomorphism between the group manifold $R=AN$ and the Riemannian 
symmetric space $M=G/K$~:
\begin{equation}\label{AN}
R\stackrel{\sim}{\to}G/K~:(a,n)\to anK.
\end{equation}
Observe that the vector space $\p\subset\g$ is naturally identified with 
the tangent space $T_{K}(M)$.
\begin{dfn}
The symmetric space $M=G/K$ is {\bf Hermitian} if there exists an 
endomorphism $J\in End(\p)$ of the vector space $\p$ such that~:
\begin{enumerate}
\item[(i)] $J^{2}=-id_{\p}$
\item[(ii)] $B(JX,JY)=B(X,Y)\quad\forall X,Y\in\p$
\item[(iii)] $ad(k)\circ J=J\circ ad(k)\quad\forall k\in\k$.
\end{enumerate}
\end{dfn}
In this case, the tensors $J$ and $B$ on $\p=T_K(M)$ globalize 
to $M$ respectively as a complex structure $\j$ and a Riemannian 
metric ${\bf g}$ on $M$. Moreover, the 2-form $\omega$ on $M$ defined 
by
$$
\omega_{x}(X,Y)={\bf g}_{x}(\j X,Y)\quad x\in M; X,Y\in T_{x}(M)
$$
is a $G$-invariant symplectic structure on $M$.

\subsection{(Co)Adjoint orbits}\label{CAO}
We now realize our Hermitian symmetric space $M=G/K$ as a coadjoint 
orbit in $\g^{\star}$ or equivalently, using the Killing form $B$,
as an adjoint orbit in $\g$.\\
Let us first denote by $\Omega\in\bigwedge^{2}(\g^{\star})$ the 
skewsymmetric 2-form on $\g$ defined by 
\begin{eqnarray*}
\Omega(X,Y)=B(JX,Y)\quad X,Y\in\p;\\
\Omega(\k,\g)=0.
\end{eqnarray*}
One observes, using the $\k$-invariance, that the 2-form $\Omega$ 
is a Chevalley 2-cocycle for the trivial representation of $\g$ on 
$\R$~:
$$
\Omega\in C^{2}_{Chevalley}(\g,\R);\qquad \delta\Omega=0.
$$
Whitehead's lemmas then  tell us that there exists an element 
$\xi_{0}\in\g^{\star}$ such that
\begin{equation}\label{CHEV}
\delta\xi_{0}=\Omega.
\end{equation}
Equivalently, one gets $Z_{0}\in\g$ such that
$$
B(Z_{0},\, . \,)=\xi_{0}.
$$
Observe that the definition of $\Omega$ implies
\begin{enumerate}
\item[(i)] $Z_{0}\in\z(\k)$,
\item[(ii)] $J=\pm ad(Z_{0})|_{\p}$ 
\end{enumerate}
where $\z(\k)$ denotes the center of $\k$.
This, together with a little more work (see \cite{Bi1, K}) provides the following ``Hamiltonian'' 
description of Hermitian symmetric spaces according to Kostant's 
classification of homogeneous Hamiltonian spaces.
\begin{prop}
\begin{enumerate}
\item[(i)] A symmetric space $G/K$ is Hermitian if and only if 
$\z(\k)\neq 0$.
\item[(ii)] In this case, $\dim(\z(\k))=1$ and the map 
$$
\z(\k)\backslash\{0\}\to\bigwedge {}^{2}(\p^{\star})~:Z\to\delta B(Z,\, . 
\,)|_{\p\times\p}
$$
induces a bijection onto the set of $\k$-invariant bilinear symplectic 
forms on $\p$.
\item[(iii)] The Hermitian symmetric space $(M,\omega,\j)$ is then 
realized as the adjoint orbit $\o=Ad(G)Z\subset\g$ of $Z$ via
$$
M=G/K\stackrel{\sim}{\to}\o~:gK\to Ad(g)Z.
$$
Under this identification, the symplectic form $\omega$ corresponds to 
the {\bf Kostant symplectic form} $\omega^{\o}$ on $\o$ defined by 
\begin{equation}\label{KOSTANT}
\omega^{\o}_{x}(X^{\star},Y^{\star})=-B(x,[X,Y])
\end{equation}
where $x\in\o\subset\g; X,Y\in\g$ and where $X^{\star}$ denotes the 
fundamental vector field associated to $X\in\g$ on $\o$~:
$$
X^{\star}_{x}=\frac{d}{dt}|_{0}Ad(\exp(-tX))x.
$$
\end{enumerate}
\end{prop}
\subsection{A class of K\"{a}hlerian groups and their Iwasawa 
coordinates}
When in a Hermitian situation, the diffeomorphism~(\ref{AN}) endows 
the group manifold $R=AN$ with the transported symplectic form $\omega$ 
coming from $M$. The symplectic form $\omega$ on $R$ is then invariant 
under the left action $$L~:R\times R\to R~:(x,y)\to xy=L_{x}y.$$
A Lie group with a left invariant symplectic structure is called a 
{\bf symplectic group} \cite{L}. 

\begin{prop}\label{IW}
Let $M=G/K$ be an irreducible Hermitian symmetric space of the non-compact 
type. Let $(R=AN,\omega)$ be the associated symplectic group via the 
isomorphism~(\ref{AN}). Denote by $\r=\a\oplus\n$ its Lie 
algebra. Then, the map
$$
\r=\a\oplus\n\stackrel{\iw}{\to} R~:(a,n)\to\exp(a)\exp(n)
$$
is a global diffeomorphism called {\bf Iwasawa coordinates}. 
Moreover, through the map $\iw$, the 
symplectic form reads
$$
(\iw^{\star}\omega)_{r}(u,v)=B\left( 
Z,\left[ Ad(exp(-n))u_{\a}+F(ad(n))u_{\n},\, 
(u \leftrightarrow v) \right]\right)
$$
where
$r=(a,n)\in\r\subset\g\mbox{ and } u,v\in T_{r}(\r)=\r\subset\g$,
where we write $u=u_{\a}+u_{\n}$ according to the decomposition 
$\r=\a\oplus\n$ and where $F$ is the analytic function defined by
$$
F(z)=\frac{1-e^{-z}}{z}\quad (z\in\C).
$$
\end{prop}
\Pf
Let $Z\in\z(\k)$ be as in Section~\ref{CAO} and denote by 
$\varphi:\r\to\o=Ad(G)Z\subset\g$ the global diffeomorphism defined 
by 
$$
\varphi(r)=Ad(\iw(r))Z.
$$
Identifying $T_{r}\r$ with $\r\subset\g$, one has for $u\in T_{r}\r$~:
\begin{eqnarray*}
\varphi_{\star_{r}}(u)=\ddto Ad(\exp(a+tu_{\a})\exp(n+tu_{n}))Z=\\
-u_{\a}^{\star}|_{\varphi(r)}+Ad(\exp a)\ddto Ad(\exp n)Ad(\exp -n)
Ad(\exp(n+tu_{\n}))Z=\\
-u_{\a}^{\star}|_{\varphi(r)}+Ad(\iw(r))\ddto Ad(CBH(-n,n+tu_{\n}))Z
\end{eqnarray*}
where $CBH$ is the Campbell-Backer-Hausdorff function for the group 
$N$ ($\exp x.\exp y=\exp CBH(x,y)$). Now, since
$$
\ddto CBH(-n,n+tu_{\n})=F(ad(n))u_{\n} \mbox{ (see \cite{H})},
$$
one gets
\begin{equation}
\varphi_{\star_{r}}(u)=-u_{\a}^{\star}|_{\varphi(r)}-
\left(Ad(\iw(r))F(ad(n))u_{\n}\right)^{\star}_{\varphi(r)}.
\end{equation}
Hence, using the $Ad$-invariance of the Killing form and the 
fact that $Ad(\exp a)|_{\a}=id_{\a}$, Formula~(\ref{KOSTANT}) yields 
the result.
\EPf

On every Hermitian symmetric space of the non-compact type $M$, there 
actually exists a global {\bf Darboux chart} i.e. a coordinate system 
where the symplectic structure reads constantly. Indeed, 
one can realize $M$ as a coadjoint orbit of $R$~: 
the orbit of the element $\iota^{\star}\xi_{0}$ where 
$\r\stackrel{\iota}{\to}\g$ is the canonical injection. A result of 
Pedersen \cite{P} then states that on the universal covering of 
every coadjoint orbit of a solvable Lie group there exists a global 
Darboux chart.\\
We will show that, at least in the rank one case, the Iwasawa 
coordinates explicitly yield such a global Darboux chart. Before 
this, 
we establish the following lemma which will be useful further on.
\begin{lem}\label{AZN}
Let $\g$ be a simple Lie algebra with Iwasawa decomposition 
$\g=\k\oplus\a\oplus\n$. Assume $\z(\k)\neq 0$. Then, 
$$
\dim \a\geq\dim\z(\n)
$$
where $\z(\n)$ denotes the center of $\n$.
\end{lem}
\Pf
Let $\r\stackrel{\iota}{\to}\g$ be the canonical inclusion and let 
$\xi_0\in\g^\star$ be such that $\delta\xi_0=\Omega$ (cf. (\ref{CHEV})).
Then, since $\k\cap\r=0$, the radical of $\delta\iota^\star\xi_0$ in 
$\r$ is trivial. Moreover, if $V$ denotes the radical of $\Omega$ in $\n$,
one has $\z(\n)\subset V$. Indeed, if
$z\in\z(\n)$, one has $\Omega(\n,z)=\xi_0[z,\n]=0$. Observe now that 
the map
$$
V\to\a^\star~:v\to\Omega(v,\, . \,)|_\a
$$
is injective. Indeed, let $v\in V$ be such that $\Omega(v,\a)=0$. Then, 
$0=\Omega(v,\a\oplus\n)=\delta\iota^\star\xi_0(v,\r)$
hence $v=0$. Thus 
$\dim\a^\star=\dim\a\geq\dim V\geq\dim\z(\n)$.
\EPf
\subsection{Rank one}
\begin{prop}\label{DARBOUX}
Let $M=G/K$ and $(R=AN,\omega)$ be as in Proposition~\ref{IW}. Assume 
$\dim\a=1$. Then, the Iwasawa coordinates $\r\stackrel{\iw}{\to}R$ define 
a global Darboux chart on $(R=AN,\omega)$ i.e. $\iw^\star(\omega)$ 
is a constant bilinear 2-form on the vector space $\r$.
\end{prop}
Before passing to the proof, we first recall the following classical result about the 
structure of $\r$.
\begin{lem}\label{STRR}
Assume $\dim\a=1$ and $\dim\g>3$. Then
\begin{enumerate}
\item[(i)] $\Phi=\{\pm\alpha,\pm2\alpha\}$;
\item[(ii)] $\n=\g_{\alpha}\oplus\g_{2\alpha}$ and $\g_{2\alpha}=\z(\n)$;
\item[(iii)] $\dim\z(\n)=\dim\a=1$;
\item[(iv)] there exists an element $E\in\z(\n)$ such that the Lie
bracket on $\n$ reads $$[x,y]=\Omega(x,y)E\quad x,y\in \n.$$
The subspaces $\a\oplus\z(\n)$ and $\g_\alpha$ are symplectic and orthogonal 
in $(\r,\Omega)$. In particular, $\n$ is a Heisenberg algebra.
\end{enumerate}
\end{lem}
\Pf
Since $\dim\a=1$, every root is a multiple of a given one, say $\alpha$. A classical 
lemma \cite{H} tells us that $\Phi\subset\{\pm\alpha,\pm2\alpha\}$. The hypothesis 
$\dim\g>3$ together with Lemma~\ref{AZN} imply (i), (ii) and (iii). 
For (iv), we  will first prove that $Z=\sigma E_0+T_0+E_0$ with 
$T_0\in\z_\k(\a)$ and $E_0\in\z(\n)$. Indeed, let $a\in\a$ and $n,n'\in\n$ 
be such that $Z=\sigma n+a+T_0+n'$ according to the 
decomposition~(\ref{BOREL}). Then $B(Z,\a)=0$ implies $a=0$; hence 
$\sigma n + n'\in\k$. Thus $n=n'$ because $\n\cap\k=0$.
Now let $z\in\z(\n)$. Then
$[Z,z+\sigma z]=0=[T_0,z]+\sigma[T_0,z]+[\sigma n,z]+\sigma[\sigma n,z]$.
Since $[\g_{0},\z(\n)]\subset\z(\n)$, one gets $[T_0,z]=0$. 
Writing $n=n_\alpha+n_{2\alpha}$ according to the decomposition 
$\n=\g_\alpha\oplus\g_{2\alpha}$, one gets $[\sigma n_\alpha,z]+
[\sigma n_{2\alpha},z]\in\p$. Hence, since $[\sigma 
n_{2\alpha},z]\in\g_0$, one has $[\sigma n_\alpha,z]\in\g_\alpha\cap\p$. 
This last intersection being zero since $[\a,\p]\subset\k$ implies 
$\n\cap\p=0$. Therefore $n=n_{2\alpha}=E_0\in\z(\n)$ and one gets the desired form 
for $Z$.\\
Now, one has 
$\Omega(\a,\g_\alpha)=B(Z,[\a,\g_\alpha])=B(\sigma E_0+E_0,\g_\alpha)=0.$
Therefore $\g_\alpha=\left(\a\oplus\z(\n)\right)^{\perp_\Omega}$ is 
symplectic and the table of $\n$ reads
$[x,y]=\frac{1}{B(\sigma E_0,E_0)}\Omega(x,y)E_0.$
\EPf

\noindent{\sl Proof of Proposition~\ref{DARBOUX}.} Assume the rank to be 
one and $\dim \g>3$. By distributing the terms of the Taylor series of function $F$ (cf. 
Proposition~\ref{IW}), one 
gets~:
$$
\left[ Ad(exp(-n))u_{\a}+F(ad(n))u_{\n},\, 
(u \leftrightarrow v) \right]= 
$$

$$
[u_\a-[n,u_\a]+{\bf c.t.}\, ,
v_\n-\frac{1}{2}[n,v_\n]]+[u_\n-\frac{1}{2}[n,u_\n]\, ,v_\a-[n,v_\a]]+
$$

$$
+[u_\n+{\bf c.t.}\, ,v_\n+{\bf c.t.}]\quad({\bf c.t.}=\mbox{central terms 
in } \n)=
$$

$$
=[u_\a,v_\n]+[u_\n,v_\a]+[u_\n,v_\n]-\frac{1}{2}[u_\a,[n,v_\n]]-
[[n,u_\a],v_\n]+[u_\n,[v_\a,n]]-\frac{1}{2}[[n,u_\n],v_\a]=
$$

$$
\mbox{(because $[n,v_\n]$ and $[n,u_\n]$ are central in $\n$)}
$$
\begin{equation}\label{EQ1}
=[u,v]-\alpha(u_\a)[n,v_\n]-[[n,u_\a],v_\n]+[u_\n,[v_\a,n]]+\alpha(v_\a)[n,u_\n].
\end{equation}
Now, observe that for $n,N\in\n$ and $A\in\a$, one has
\begin{eqnarray*}
[[n,A],N]=[[n_\alpha+n_{2\alpha},A],N]=
[N,\alpha(A)n_\alpha+2\alpha(A)n_{2\alpha}]=\\
\alpha(A)[N,n_\alpha]=\alpha(A)[N,n]\quad(\mbox{because } 
\g_{2\alpha}=\z(\n)).
\end{eqnarray*}
Hence (\ref{EQ1}) becomes
$[u,v]-\alpha(u_\a)[n,v_\n]-\alpha(u_\a)[v_\n,n]+\alpha(v_\a)[u_\n,n]
+\alpha(v_\a)[n,u_\n]=[u,v]$. The case $\dim\g=3$ is similar and simpler.
\EPf

\section{Guessing the product formula}\label{Scarlatti}
Star products have been introduced in \cite{BFFLS} as an autonomous 
formulation of Quantum Mechanics. 
\begin{dfn}
A {\bf star product} on a symplectic manifold $(M,\omega)$ is an associative $\C$-bilinear 
multiplication $\star_\nu$ on the space of formal power series $C^\infty(M)[[\nu]]$ such that, 
for all $u,v\in C^\infty(M)$, one has
\begin{enumerate}
\item[(i)] 
$$
u\star_\nu v=\sum_{k=0}^\infty c_k(u,v)\nu^k
$$ where the $c_k$'s are bidifferential operators on $C^\infty(M)$;
\item[(ii)] $c_0(u,v)=uv$;
\item[(iii)] $c_1(u,v)=\frac{1}{2}\{u,v\}$ where $\{\, , \, \}$ is the Poisson bracket 
associated to the symplectic form $\omega$;
\item[(iv)] $u\star_\nu1=1\star_\nu u=u$.
\end{enumerate}
\end{dfn}
\begin{ex}
{\rm
Let $(V,\Omega)$ be a symplectic vector space of dimension $2n$ and write 
$\Omega(x,y)=<x,Jy>$ with $J\in so(V)$. Then evaluating Rieffel's product 
(\ref{RIEFFEL})
of two compactly supported functions $u,v\in C^\infty_c(V)$ at a point $x\in V$, 
one re-finds the old expression of the {\bf Weyl product}~:
\begin{equation}\label{WEYL}
(u\star^{W}_{\hbar}v)(x)\stackrel{\mbox{def.}}{=}
(u\star_{\hbar J}v)(x)=\hbar^{-2n}\int_{V\times V} e^{\frac{2i}{\hbar}S^0(x,y,z)}u(x)\, v(y)\, dy\, dz
\end{equation}
with 
$$
S^0(x,y,z)=<x,Jy>+<y,Jz>+<z,Jx>,
$$
as explained in the introduction. Recall that the Schwartz space $\s(V)$ is 
stable under Weyl's product \cite{Ha}. Performing a stationary phase method on the 
oscillatory integral (\ref{WEYL}), 
one gets an asymptotic expansion in the parameter $\nu=\frac{\hbar}{2i}$~:
\begin{equation}\label{MOYAL}
u\star_{\hbar J}v\sim
 \, uv+\nu\{u,v\}+\sum_{k=2}^{\infty}\frac{\nu^{k}}{k!}\sum_{\begin{array}{c} 
i_{1}\ldots i_{k}\\ j_{1}\ldots j_{k}
\end{array}} \Omega^{i_{1}j_{1}}\ldots \Omega^{i_{k}j_{k}} 
\partial_{i_{1}\ldots i_{k}}u.\partial_{j_{1}\ldots j_{k}}v.
\end{equation}
The expression in the RHS of (\ref{MOYAL}) actually defines a star product on 
$C^\infty(M)[[\nu]]$ called the {\bf Moyal star product}. It will be denoted by $\star_\nu^M$.
}
\end{ex}
\subsection{$AN$-covariant Moyal star products}\label{Vivaldi}
In this section, we adapt to our situation old techniques from star
representation theory \cite{AC, BFFLS, Fr}. First, we recall the notion of
covariant star product.
\begin{dfn}\label{HAM}
let $(M,\omega)$ be a symplectic manifold on which a connected Lie
group $G$ acts in a strongly Hamiltonian way. Let $\g$ be the Lie
algebra of $G$ and denote by 
$$
\g\stackrel{\lambda}{\longrightarrow} C^\infty(M): X\to \lambda_X
$$
the (dual) {\bf moment map} i.e. 
$$
i_{X^\star}\omega
=-d\lambda_X.
$$
A star product $\star_\nu$ on $C^\infty(M)[[\nu]]$ is said to be
{\bf $G$-covariant} if 
$$
[\lambda_X,\lambda_Y]_{\star_\nu}\stackrel{\mbox{def.}}{=}
\lambda_X\star_\nu\lambda_Y-\lambda_Y\star_\nu\lambda_X=2\nu\{\lambda_X,\lambda_Y\}.
$$
\end{dfn}
\begin{prop}\label{MOMENT}
Within the assumptions and notations 
of Proposition~\ref{IW}, let $R\times\r\stackrel{\tau}{\to}\r$ be
the action defined by
$$
\tau_g(r)=\iw^{-1}g\,\iw(r).
$$
Assume $\dim\r\geq4$ (i.e. $\dim\g>3$). 
Then, this action is Hamiltonian with respect to the constant
symplectic structure $\Omega$ on $\r$. Moreover, the Hamiltonian
functions associated to the infinitesimal action are
$$
\begin{array}{cccc}
\lambda_A(r) & = & 2\alpha(A)B(\sigma E,E)n_E & (A\in\a),\\
\lambda_y(r) & = & e^{-\alpha(a)}\Omega(n,y) & (y\in\g_\alpha),\\
\lambda_E(r) & = & e^{-2\alpha(a)}B(\sigma E,E), & 
\end{array}
$$
where $r=(a,n)$ and $n=n_\alpha+n_EE$ according to the decomposition
$\n=\g_\alpha\oplus\R E$ (cf. Lemma~\ref{STRR}). In particular, every such
Hamiltonian is linear in $\n$, therefore the Moyal star product on
$(\r,\Omega)$ is $R$-covariant.
\end{prop}
\Pf
\begin{eqnarray*}
\lambda_A(r) & = & B(Ad(\exp a\exp n)Z_0,A)\\
	     & = & B(Z_0,Ad(\exp -n)A)\\
	     & = & B(Z_0, A-[n,A]+\frac{1}{2}[n,[n,A]])\\
	     & = & B(\sigma E, A-[n,A]+\frac{1}{2}[n,[n,A]])\\
	     & = & B(\sigma E, 2\alpha(A)n_E
	     E+\frac{1}{2}[\alpha(A)n_\alpha,n_\alpha])\\
             & = & 2\alpha(A)B(\sigma E,E)n_E;\\
\lambda_y(r) & = &  B(Ad(\exp a\exp -n)Z_0,y)\\
	     & = & e^{-\alpha(a)}B(Z_0,Ad(\exp -n)y)\\
	     & = & e^{-\alpha(a)}B(Z_0,y-[n,y])\\
	     & = &
	     e^{-\alpha(a)}B(Z_0,-[n,y])=e^{-\alpha(a)}\Omega(n,y);\\
\lambda_E(r) & = &  B(Ad(\exp a\exp n)Z_0,E)\\
	     & = & e^{-2\alpha(a)}B(Z_0,Ad(\exp -n)E)\\
	     & = & e^{-2\alpha(a)}B(Z_0,E).
\end{eqnarray*}

\EPf

When covariant, a star product yields a representation of $\g$ on
$C^\infty(M)[[\nu]]$~:
$$
\g\stackrel{\rho_\nu}{\longrightarrow}\mbox{End}(C^\infty(M)[[\nu]])
$$

$$
\rho_\nu(X)u=\frac{1}{2\nu}[\lambda_X,u]_{\star_\nu}.
$$
In order to compute the representation $\rho_\nu$ in our context, we
observe
\begin{lem}
let $(\r,\Omega)$ be a symplectic vector space. Let $\u$ be a
codimension 2 symplectic subspace of $\r$ and let $A$ and $E$ be
generators of $\u^\perp$. Then, for every linear from $\mu\in\u^\star$ 
and every smooth function $\epsilon\in C^\infty(\R.A)$, one has
$$
\Omega^{i_1j_1}...\Omega^{i_kj_k}\partial_{i_1...i_k}(\epsilon\otimes\mu)\partial_{j_1...j_k}(u)
=k\,\partial_A^{k-1}\epsilon\,\partial_{{}^\sharp\mu}\partial_E^{k-1}u+\mu\,
\partial_A^{k}\epsilon\,\partial_E^{k}u
\quad (u\in C^\infty(\r)),
$$
where ${}^\sharp\mu$ is defined by $\Omega({}^\sharp\mu,\, . \,)=-\mu$.
\end{lem}
\Pf
One has $\partial_E(\epsilon\otimes\mu)=0$ and
$\partial^\ell_\u(\epsilon\otimes\mu)=0$ as soon as $\ell\geq2$. Hence
$\partial_{i_1...i_k}(\epsilon\otimes\mu)\neq 0$ only if the $k$-tuple
$(i_1...i_k)$ 
contains either one or zero element of $\u$; all the other ones being
$A$'s. There are $k$ such $k$-tuples for a given element of
$\u$. Therefore, the only $(j_1...j_k)$'s yielding non zero
contributions in the LHS contain either one or zero element of $\u$
(conjugated with the one in the corresponding $(i_1...i_k)$) and
$E$'s. Therefore, one gets
$$
LHS=k\Omega^{\alpha\beta}_\u\,\partial_A^{k-1}\partial_{\alpha}(\epsilon\otimes\mu)\,
\partial_\beta\,\partial_E^{k-1}u+\partial_A^{k}(\epsilon\otimes\mu)\,\partial_E^{k}u.
$$
One concludes using $\Omega^{\alpha\beta}\partial_\alpha\mu\,\partial_\beta=\partial_{{}^\sharp\mu}$.

\EPf

This implies that for $y\in\g_\alpha$, one has 
\begin{eqnarray*}
\frac{1}{2\nu}[\lambda_y,u]_{\star^M_\nu} & = & \frac{2}{2\nu}\sum_k\frac{\nu^{2k+1}}{(2k+1)!}
\left\{
(2k+1)(-\alpha(A))^{2k}e^{-\alpha(a)}\partial_y\partial_E^{2k}u +
\Omega(n,y)(-\alpha(A))^{2k+1}e^{-\alpha(a)}\partial_E^{2k+1}u
\right\}\\
  & = & e^{-\alpha(a)}\cosh\left(\nu(-\alpha(A))\partial_E\right)\partial_yu+
\frac{1}{\nu}\Omega(n,y)\,\sinh\left(\nu(-\alpha(A))\partial_E\right)e^{-\alpha(a)}u.
\end{eqnarray*}
Also
$$
\frac{1}{2\nu}[\lambda_A,u]_{\star^M_\nu}=\partial_Au
$$
and 
\begin{eqnarray*}
\frac{1}{2\nu}[\lambda_E,u]_{\star^M_\nu} & = & \frac{2}{2\nu}\sum_k\frac{\nu^{2k+1}}{(2k+1)!}
B(\sigma E,E)(-2\alpha(A))^{2k+1}e^{-2\alpha(a)}\partial_E^{2k+1}u\\
  & = & \frac{1}{\nu}B(\sigma E,E)\sinh\left(\nu(-2\alpha(A))\partial_E\right)e^{-2\alpha(a)}u.
\end{eqnarray*}
Regarding these expressions, it is tempting to take the partial
Fourier transform in the $E$-variable in order to obtain a so called
``multiplicative representation".\\

\noindent Writing an element $r\in\r$ as 
$$
r=aA+x+zE\mbox{ with }x\in\g_\alpha,
$$
we set, for (reasonable) $u\in C^\infty(\r)$,
$$
F(u)(a,x,\xi)=\hat{u}(a,x,\xi)=\int_{\z(\n)}e^{-i\xi z}u(aA+x+zE)\,dz.
$$
One then has $F(\partial_Eu)=i\xi\hat{u}$ which yields
\begin{eqnarray*}
F(\rho_\nu(y)u)=e^{-\alpha(a)}\cosh\left(\nu(-\alpha(A))i\xi\right)\partial_y\hat{u}+
\frac{1}{\nu}\Omega(x,y)\,\sinh\left(\nu(-\alpha(A))i\xi\right)e^{-\alpha(a)}\hat{u};\\
F(\rho_\nu(A)u)=\partial_A\hat{u};\\
F(\rho_\nu(E)u)=\frac{1}{\nu}B(\sigma E,E)\sinh\left(\nu(-2\alpha(A))i
\xi\right)e^{-2\alpha(a)}\hat{u}.
\end{eqnarray*}
Choosing $A$ such that $\alpha(A)=1$ and setting
$\nu=\frac{\hbar}{2i}$, one gets
\begin{eqnarray*}
\hat{\rho}_\hbar(y)\hat{u}=e^{-a}\cosh\left(\frac{\hbar}{2}\xi\right)\partial_y\hat{u}+
\frac{2i}{\hbar}\Omega(x,y)\,\sinh\left(\frac{\hbar}{2}\xi\right)e^{-a}\hat{u};\\
\hat{\rho}_\hbar(A)\hat{u}=\partial_A\hat{u};\\
\hat{\rho}_\hbar(E)\hat{u}=-\frac{2i}{\hbar}B(\sigma
E,E)\sinh\left(\hbar\xi\right)
e^{-2a}\hat{u},
\end{eqnarray*}
where $\hat{\rho}_\hbar$ is the representation of $\g$ defined by 
$$
\hat{\rho}_\hbar(X)\hat{u}=F(\rho_{\frac{\hbar}{2i}}u).
$$
This is a multiplicative representation. We now change the coordinates
following 
\begin{equation}
(a,x,\xi)=(a',\cosh(\frac{\hbar}{2}\xi')x',\xi')\stackrel{\mbox{def.}}{=}\varphi_\hbar(a',x',\xi').
\end{equation}
This yields
\begin{eqnarray*}
\hat{\rho}_\hbar(y)f(a',x',\xi')=e^{-a'}\partial'_yf-\Omega(x',y)e^{-a'}
\frac{i}{\hbar}\sinh(\hbar\xi')f;
\end{eqnarray*}
the rest being unchanged.
\begin{dfn}
We denote by $c_\hbar\in\Omega^1(\r)$ the smooth one-form on $\r$
defined by
$$
(c_\hbar)_{(a,x,\xi)}(X)=-e^{-a}\frac{\sinh(\hbar\xi)}{\hbar}\left( 
\Omega(x,X)+2e^{-a}B(Z_0,X) 
\right)\quad X\in T_{(a,x,\xi)}(\r).
$$
\end{dfn}
One then gets
\begin{lem}
Under the transformation $\varphi^\star_\hbar\circ F$, the star
representation of $\r=Lie(R)$ on $C^\infty(\r)[[\hbar]]$ is
multiplicative and reads as
$$
\pi_\hbar(X)f(a,x,\xi)\stackrel{\mbox{def.}}{=}\left(e^{-a}X_\alpha+X_\a\right).f
+ ic_\hbar(X)f\qquad (f\in C^{\infty}(\r))
$$
where $X=X_\a+X_\alpha+X_E$ according to the decomposition
$\r=\a\oplus\g_\alpha\oplus\R E$ and where
$$
\pi_\hbar(X)=\varphi^\star_\hbar\circ\hat{\rho}_\hbar(X)\circ(\varphi_\hbar^{-1})^\star
=\varphi^\star_\hbar\circ F\circ \rho_{\frac{\hbar}{2i}}(X)\circ F^{-1}\circ(
\varphi_\hbar^{-1})^\star.
$$
\end{lem}

\subsection{The $\z_\hbar$-transform}\label{Haendel}
\begin{dfn} For $u\in C^{\infty}(\r)$ integrable, we define the {\bf $\z_\hbar$-transform} 
by
$$
\left(\z_\hbar(u)\right)(a,x,\xi)=
\int e^{-\frac{i}{\hbar}\sinh(\hbar\xi)z}u(a,x,z) \, dz.
$$
The formal (commutative) product obtained by transporting the pointwise multiplication 
of functions via $\z_\hbar$ is denoted by $\bullet_\hbar$~:
$$
f\bullet_\hbar g\stackrel{\mbox{def.}}{=} \z_\hbar(\z_\hbar^{-1}f.
\z_\hbar^{-1}g)
$$
(whenever this expression makes sense).
\end{dfn}

\begin{thm}
Under representation $\pi_\hbar$ , the Lie algebra $\r$ acts by derivations 
with respect to the commutative product $\bullet_\hbar$ i.e. one 
has formally
$$
\pi_\hbar(f\bullet_\hbar g)=(\pi_\hbar(X)f)\bullet_\hbar g+f\bullet_\hbar(\pi_\hbar(X)g).
$$
\end{thm}

\Pf
It is sufficient to prove that, for all
$X\in\r$,
\begin{equation}\label{VF}
X^\hbar\stackrel{\mbox{def.}}{=}\z_\hbar^{-1}\circ
m_{c_\hbar(X)}\circ \z_\hbar\mbox{ is a vector field on }\r
\end{equation}
($m_{c_\hbar(X)}$ denotes the multiplication by $c_\hbar(X)\, : \, 
m_{c_\hbar(X)}(f)=c_\hbar(X)f$). Indeed, if (\ref{VF}) holds one has
\begin{eqnarray*}
c_\hbar(X)f\bullet_\hbar g+f\bullet_\hbar c_\hbar(X)g=\\
\z_\hbar(X^\hbar\z_\hbar^{-1}f.\z_\hbar^{-1}g
+\z_\hbar^{-1}f.X^\hbar\z_\hbar^{-1}g)=\\
\z_\hbar X^\hbar(\z_\hbar^{-1}f.
\z_\hbar^{-1}g)=\\
c_\hbar(X)(f\bullet_\hbar g).
\end{eqnarray*}
Therefore, since the vector part $\tilde{X}=e^{-a}X_\alpha+X_\a$ 
of $\pi_\hbar(X)$ does not involve the $E$-variable, its action commutes 
with the $\z_\hbar$-transform and one gets
\begin{eqnarray*}
\pi_\hbar(f\bullet_\hbar g)=\\
\z_\hbar(\tilde{X}(\z_\hbar^{-1}f.\z_\hbar^{-1}g))+ic_\hbar(X)(f\bullet_\hbar 
g)=\\
\z_\hbar((\z_\hbar^{-1}\tilde{X}f)\z_\hbar^{-1}g
+\z_\hbar^{-1}f(\z_\hbar^{-1}\tilde{X}g))+i(c_\hbar(X)f)\bullet_\hbar 
g+if\bullet_\hbar(c_\hbar(X)f)=\\
(\tilde{X}f)\bullet_\hbar g+f\bullet_\hbar(\tilde{X}g)+i(c_\hbar(X)f)\bullet_\hbar 
g+if\bullet_\hbar(c_\hbar(X)f)=\\
(\pi_\hbar(X)f)\bullet_\hbar g+f\bullet_\hbar(\pi_\hbar(X)g).
\end{eqnarray*}
We now prove assertion (\ref{VF}). For $y\in\g_\alpha$, one has
$$
c_\hbar(y)\z_\hbar(u)=e^{-a}\Omega(y,x)\frac{\sinh(\hbar\xi)}{\hbar}\int 
e^{-\frac{i}{\hbar}\sinh(\hbar\xi)z}u(a,x,z) \, 
dz=-ie^{-a}\Omega(y,x)\z_\hbar(\partial_zu).
$$
Hence, since the $\z_\hbar$-transform only involves the $E$-variable,
$$
\z_\hbar^{-1}\circ
m_{c_\hbar(y)}\circ \z_\hbar=-ie^{-a}\Omega(y,x)\partial_z.
$$
For $A\in\a$, one has
$c_\hbar(A)=0$. For $E\in\z(\n)$, one has
$$
\z_\hbar(\partial_Eu)=\int e^{-\frac{i}{\hbar}\sinh(\hbar\xi)z}\partial_Eu \, dz
=\frac{i}{\hbar}\sinh(\hbar\xi)\z_\hbar 
(u)=-i\frac{e^{2a}c_\hbar(E)}{2B(Z_0,E)}\z_\hbar 
(u),
$$
hence
$$
\z_\hbar^{-1}c_\hbar(E)\z_\hbar(u)=2iB(Z_0,E)e^{-2a}\partial_Eu.
$$
\EPf

If one interprets the commutative product $\bullet_{\hbar}$ as the 
underlying product to the algebra of functions on a commutative 
$\hbar$-dependent manifold, say $M_{\hbar}$, its invariance under 
$\hat{\rho}$ tells us that $\g$ is realized via $\hat{\rho}$  as a 
subalgebra of tangent vector fields over $M_{\hbar}$. \\
At this level, we want
\begin{enumerate}
\item[(a)] to identify the infinitesimal action $\g\to\Gamma(TM_\hbar)$;
\item[(b)] to identify the product on $C^\infty(\r)[[\hbar]]$ defined by 
$$
u\star_\hbar v=T_\hbar^{-1}(T_\hbar u\star^M_{\frac{\hbar}{2i}}T_\hbar v)
$$
with $T_\hbar=F^{-1}\circ(\varphi^{-1})^\star\circ\z_\hbar$.
\end{enumerate}

Formally $\star_\hbar$ is indeed a quantization of $(\r,\Omega)$ since 
$\lim_{\hbar\to 0}T_\hbar=id$. The following proposition answers question (a).

\begin{prop}
For all $X\in\r$, one has
$$
\z_\hbar^{-1}\circ\pi_\hbar(X)\circ\z_\hbar=X^\star,
$$
whenever this expression makes sense.
\end{prop}

\Pf
First, one has $$\z_\hbar^{-1}\pi_\hbar(A)\z_\hbar(u)=
\z_\hbar^{-1}\partial_A\z_\hbar(u)=\partial_Au=A^\star u.$$ Second, for $X=y+xE\in\n\quad(y\in\g_\alpha)$, one has
$$
X^\star u=\{\lambda_X,u\}=e^{-a}\partial_yu-\Omega(x,y)e^{-a}\partial_Eu-
2B(Z_0,Z)e^{-2a}\partial_Eu
$$
(cf. Proposition~\ref{MOMENT}). Hence
\begin{eqnarray*}
\z_\hbar X^\star u=\\
\int e^{-\frac{i}{\hbar}\sinh(\hbar\xi)z} \left\{ 
e^{-a}\partial_yu-\Omega(x,y)e^{-a}\partial_Eu-
2B(Z_0,Z)e^{-2a}\partial_Eu
\right\}=\\
-\frac{i}{\hbar}\sinh(\hbar\xi)e^{-a}(\Omega(x,y)+2B(Z_0,Z)e^{-a})\z_\hbar(u)
+e^{-a}\partial_y\z_\hbar(u)=\\
ic_\hbar(X)\z_\hbar(u)+e^{-a}\partial_y\z_\hbar(u)=\\
\pi_\hbar(X)\z_\hbar(u).
\end{eqnarray*}
\EPf

Therefore $M_\hbar$ can be $R$-equivariantly identified with $\r$, which 
implies that the star product $\star_\hbar$ on $(\r,\Omega)$ described 
in (b) is $R$-invariant.\\
In order to define function algebras which will be stable under the 
product $\star_\hbar$, we will transport the structure of the 
Schwartz space--- which is stable under the Weyl product (\ref{WEYL})--- via 
the ``equivalence" $T_\hbar$.
\section{Strict Quantization}\label{Bach}
In this section, we adopt the following notation. If $V$ is a finite dimensional vector space, 
we denote by $\s(V)$ (resp. $\s'(V)$) the space of Schwartz functions (resp. 
tempered distributions) on $V$.
\begin{lem}\label{PHI}
Let $\phi_\hbar~:\r\to\r$ be the diffeomorphism defined by 
$$
\phi_\hbar(a,x,\xi)=(a,\frac{1}{\cosh(\frac{\hbar}{2}\xi)}x,\frac{1}{\hbar}\sinh(\hbar\xi)).
$$
Then, one has
\begin{enumerate}
\item[(i)] $\phi_\hbar^\star\s(\r)\subset\s(\r)$,
\item[(ii)] $(\phi_\hbar^{-1})^\star\s(\r)\subset\s'(\r)$.
\end{enumerate}
\end{lem}

\Pf
For the sake of simplicity, we will only prove that, if $\phi~:\R^2\to\R^2$ 
is defined by $\phi(x,y)=(\sech(\frac{y}{2})x,\sinh(y))$, then 
$\phi^\star\s(\R^2)\subset\s(\R^2)$ and 
$(\phi^{-1})^\star\s(\R^2)\subset\s'(\R^2)$. The proofs of items (i)
and (ii) being entirely similar.\\
First, one has 
$$
\phi^{-1}(x,y)=(\frac{\sqrt{2}}{2}(1+\sqrt{1+y^2})^{\frac{1}{2}}x,\arcsinh(y)),
$$

$$
\phi_{\star_{(x,y)}}=
\left(
\begin{array}{cc}
\sech(\frac{y}{2}) & -\frac{x}{2}\tanh(\frac{y}{2})\sech(\frac{y}{2})\\
0 & \cosh(y)
\end{array}
\right) \mbox{ and}
$$

\begin{eqnarray}\label{GR}
\frac{\sqrt{2}}{2}(1+\sqrt{1+y^2})^{\frac{1}{2}}=\cosh(\frac{\arcsinh(y)}{2}).
\end{eqnarray}
Therefore, setting $p_{n,m}(x,y)=x^ny^m$, $p_{n,m}\circ\phi^{-1}$ has 
still polynomial growth. This implies that for all $n,m$,  
$\sup_a\{| p_{n,m}(a)\phi^\star u(a)|\}<\sup_a\{| p_{N,M}(a)u(a)|\}$ for 
some $N,M$. This last expression being finite if $u\in\s(\R^2)$. For 
derivatives of $\phi^\star(u)$, one needs to control the asymptotic 
behavior of $\phi_{\star_{\phi^{-1}(x,y)}}$. Formulas (\ref{GR}) imply 
that $||\phi_{\star_{\phi^{-1}(x,y)}}||$ has polynomial growth. 
This shows that $\sup\{| p_{n,m}D\phi^\star u|\}<\infty$. An induction 
argument then yields $\phi^\star\s(\R^2)\subset\s(\R^2)$.\\
Now, one wants to find $N,M\geq 0$ such that 
$$
\int_U\left| x^{-N}y^{-M}(\phi^{-1})^\star u(x,y)\right|dx\, dy
$$
is finite ($U$ is the complementary subset of some compact neighborhood of the 
origin). Changing the variables following $a\rightarrow\phi(a)$, this integral 
becomes
\begin{eqnarray*}
\int_{U'}|\frac{1}{p_{N,M}(\phi(a))}u(a)||\mbox{Jac}_\phi(a)|da=\\
\int_{U'}\frac{|2\sinh(\frac{y}{2})|}{\left|\left(\frac{x}{\cosh(\frac{y}{2})}\right)^N
(\sinh(y))^M
\right|}|u(a)|da=\\
\int_{U'}\frac{2^{1-M}\left|\sinh(\frac{y}{2})^{1-M}\right|}{x^N\cosh(\frac{y}{2})^{N-M}}
|u(a)|da
\end{eqnarray*}
which is finite as soon as $N\geq M\geq 1$ ($U'$ is of the same type as $U$). For derivatives, one needs to 
control $||\phi^{-1}_{\star_{\phi(a)}}(A)||$ i.e. the norm of the inverse 
matrix $[\phi_\star]^{-1}$, which, by looking at formulas (\ref{GR}), has 
polynomial growth.
An induction argument then yields 
$(\phi^{-1})^\star\s(\R^2)\subset\s'(\R^2)$.
\EPf

\noindent Lemma~\ref{PHI} allows us to define the following linear injection~:
$$
\begin{array}{c}
\tau_\hbar~: \s(\r)\to\s'(\r)\\
\tau_\hbar\stackrel{\mbox{def.}}{=}F^{-1} \circ(\phi^{-1}_\hbar)^\star
\circ F
\end{array}
$$
where one extends the Fourier transform to the tempered distributions. We 
then set
$$
\EE\stackrel{\mbox{def.}}{=}\tau_\hbar\s(\r)\subset\s'(\r).
$$

\begin{lem}
\begin{enumerate}
\item[(i)] $\s(\r)\subset\EE$.
\item[(ii)] The map $T_\hbar:\s(\r)\to\s(\r):\quad T_\hbar=F^{-1}\circ\phi_\hbar^\star\circ F$ 
extends to $\EE$ as a linear 
isomorphism $T_\hbar:\EE\to\s(\r)$.
\item[(iii)] One has $T_\hbar\circ\tau_\hbar=id_{\s(\r)}$ and $\tau_\hbar\circ 
T_\hbar|_{\s(\r)}=id_{\s(\r)}$.
\end{enumerate}
\end{lem}

\Pf
For $u\in\s(\r)$, one has $T_\hbar(u)\in\s(\r)$ hence $\tau_\hbar 
T_\hbar(u)=u\in\EE$. The rest is obvious.
\EPf

\begin{thm}\label{THM}
\begin{enumerate}
\item[(i)] For $a,b\in\EE$, the formula
$$
a\star_\hbar b \stackrel{\mbox{def.}}{=}\tau_\hbar(T_\hbar a \star^W_\hbar 
T_\hbar b)
$$
defines an associative algebra structure on $\EE$ ($\star^W_\hbar$ denotes 
the Weyl product on $\s(\r)$, see Formula~(\ref{WEYL})).
\item[(ii)] For $u,v\in\s(\r)\subset\EE$, the product $\star_\hbar$ reads
$$
(u\star_\hbar v)(a_0,x_0,z_0)=
$$

$$
\frac{1}{\hbar^{2\dim\r}}\int_{\r\times\r}
\cosh(2(a_1-a_2))\cosh(a_2-a_0)\cosh(a_0-a_1)\times
$$

$$
\times
\exp\left(
\frac{2i}{\hbar}\left\{
S^0(\cosh(a_1-a_2)x_0, \cosh(a_2-a_0)x_1, \cosh(a_0-a_1)x_2)-
\frac{1}{2}\oint_{0,1,2}\sinh(2(a_0-a_1))z_2
\right\}
\right)\times
$$

\begin{equation}\label{PRODUCT}
\times u(a_1,x_1,z_1)\,v(a_2,x_2,z_2)\, da_1da_2dx_1dx_2dz_1dz_2
\end{equation}
where $S^0$ is the phase for the Weyl product (cf. Formula~(\ref{WEYL})) and 
where $\oint_{0,1,2}$ stands for cyclic summation.
\item[(iii)] In the Iwasawa coordinates $\r\stackrel{\iw}{\to}R$, 
the group multiplication law reads
$$
L_{(a,x,z)}(a',x',z')=\left(
a+a',e^{-a'}x+x',e^{-2a'}z+z'+\frac{1}{2}\Omega(x,x')e^{-a'}
\right).
$$
Both phase and amplitude occurring in formula (\ref{PRODUCT}) are invariant 
under the left action $L:R\times R\to R$.
\end{enumerate}
\end{thm}

\Pf
We perform the computation which leads to formula (\ref{PRODUCT}).
On the one hand, we have
$$
(T_\hbar u \star^W_\hbar 
T_\hbar v)(a_0,x_0,z_0)=\int (T_\hbar u)(a_1,x_1,z_1)\, (T_\hbar v)(a_2,x_2,z_2)
\exp\left(\frac{2i}{\hbar}
S^0(p_0,p_1,p_2)
\right)dp_1dp_2
$$
with $p_i=(a_i,x_i,z_i)$, that is 
$$
\int e^{i\Omega(\xi_1,z_1)}(\phi_\hbar^\star\hat{u})(a_1,x_1,\xi_1)
e^{i\Omega(\xi_2,z_2)}(\phi_\hbar^\star\hat{v})(a_2,x_2,\xi_2)\times
$$

$$
\begin{array}{c}
\times
\exp\left\{\frac{2i}{\hbar}\left(
\Omega(a_0,z_1)-\Omega(a_1,z_0)+\Omega(x_0,x_1)+\right.\right. \\
\left.\left.
+\Omega(a_1,z_2)-\Omega(a_2,z_1)+\Omega(x_1,x_2)
+\Omega(a_2,z_0)-\Omega(a_0,z_2)+\Omega(x_2,x_0)
\right)\right\}
\end{array}
$$

(we omit the $dp_i$'s and other such differentials)
$$
=\int\exp\left(
i\left\{ 
\Omega(\xi_1+\frac{2}{\hbar}(a_0-a_2),z_1)
+\Omega(\xi_2+\frac{2}{\hbar}(a_1-a_0),z_2)
\right\}+ \frac{2i}{\hbar}\Omega(a_2-a_1,z_0)
+\frac{2i}{\hbar}S^0(x_0,x_1,x_2)
\right)\times
$$

$$
\times (\phi_\hbar^\star\hat{u})(a_1,x_1,\xi_1) (\phi_\hbar^\star\hat{v})(a_2,x_2,\xi_2)
$$

$$
=\int\exp\left(
\frac{2i}{\hbar}\left\{ S^0(x_0,x_1,x_2)+\Omega(a_2-a_1,z_0)\right\}
\right)
(\phi_\hbar^\star\hat{u})(a_1,x_1,\frac{2}{\hbar}(a_2-a_0))
(\phi_\hbar^\star\hat{v})(a_2,x_2,\frac{2}{\hbar}(a_0-a_1)).
$$
One the second hand, one has
$$
\tau_\hbar u(a_0,x_0,z_0)=\int 
e^{i\Omega(\xi,z_0)}(\phi^{-1}_\hbar)^\star\hat{u}(a_0,x_0,\xi)\,d\xi
$$

$$
=\int 
e^{i\Omega(\xi,z_0)}\hat{u}\left(
a_0,\cosh(\frac{1}{2}\arcsinh(\hbar\xi))x_0,\frac{1}{\hbar}\arcsinh(\hbar\xi)
\right)
$$

$$
=\int \exp\left( i\left\{
\Omega(\xi,z_0)-\Omega(\frac{1}{\hbar}\arcsinh(\hbar\xi),z)
\right\}\right)
u\left(
a_0,\cosh(\frac{1}{2}\arcsinh(\hbar\xi))x_0,z
\right)\, dz\, d\xi.
$$
Therefore, one gets
$$
\tau_\hbar(T_\hbar u \star^W_\hbar 
T_\hbar v)(a_0,x_0,z_0)=\int \exp\left(
\frac{2i}{\hbar}\left\{ 
S^0(\cosh(\frac{1}{2}\arcsinh(\hbar\xi))x_0,x_1,x_2)+\Omega(a_2-a_1,z)
\right\}
\right)\times
$$

$$
\times
(\phi_\hbar^\star\hat{u})(a_1,x_1,\frac{2}{\hbar}(a_2-a_0))
(\phi_\hbar^\star\hat{v})(a_2,x_2,\frac{2}{\hbar}(a_0-a_1))
\exp\left( i\left\{
\Omega(\xi,z_0)-\Omega(\frac{1}{\hbar}\arcsinh(\hbar\xi),z)
\right\}\right)
$$
which is, changing the variables following $\eta=\arcsinh(\hbar\xi)$~:
$$
\int\exp\left(
\frac{2i}{\hbar}\left\{ 
S^0(\cosh(\frac{1}{2}\eta)x_0,x_1,x_2)+\Omega(a_2-a_1,z)
\right\}
\right)
(\phi_\hbar^\star\hat{u})(a_1,x_1,\frac{2}{\hbar}(a_2-a_0))
(\phi_\hbar^\star\hat{v})(a_2,x_2,\frac{2}{\hbar}(a_0-a_1))
\times
$$

$$
\times
\exp\left(i\left\{
\Omega(\frac{1}{\hbar}\sinh(\eta),z_0)-\Omega(\frac{1}{\hbar}\eta,z)
\right\}
\right)
\frac{1}{\hbar}\cosh(\eta)
$$

$$
=\int\exp\left(
\frac{2i}{\hbar}\left\{ 
S^0(\cosh(\frac{1}{2}\eta)x_0,x_1,x_2)+\Omega(a_2-a_1,z)
\right\}
\right)
\times
$$

$$
\times
\exp\left(
-i\left\{
\Omega(\frac{1}{\hbar}\sinh(2(a_2-a_0)),z_1)+\Omega(\frac{1}{\hbar}\sinh(2(a_0-a_1)),z_2)
\right\}
\right)
\times
$$

$$
\times
u(a_1,\cosh(a_2-a_0)^{-1}x_1,z_1)\,v(a_2,\cosh(a_0-a_1)^{-1}x_2,z_2)
\exp\left(i\left\{
\Omega(\frac{1}{\hbar}\sinh(\eta),z_0)-\Omega(\frac{1}{\hbar}\eta,z)
\right\}
\right)
\frac{1}{\hbar}\cosh(\eta)
$$
which yields, changing the variables following 
$\eta\leftarrow\frac{1}{\hbar}\eta$~:
$$
=\int\exp\left(
\frac{2i}{\hbar}
S^0(\cosh(\frac{\hbar}{2}\eta)x_0,x_1,x_2)
\right)
u(a_1,\cosh(a_2-a_0)^{-1}x_1,z_1)\,v(a_2,\cosh(a_0-a_1)^{-1}x_2,z_2)
\times
$$

$$
\times
\exp\left(
i\left\{
-\Omega(\frac{1}{\hbar}\sinh(2(a_2-a_0)),z_1)-\Omega(\frac{1}{\hbar}\sinh(2(a_0-a_1)),z_2)
+\Omega(\frac{1}{\hbar}\sinh(\hbar\eta),z_0)
\right\}
\right)
\times
$$

$$
\times
\cosh(\hbar\eta)\exp(i\Omega(\frac{2}{\hbar}(a_2-a_1)-\eta,z))
$$

$$
=\int\exp\left(
\frac{2i}{\hbar}
S^0(\cosh(a_1-a_2)x_0,x_1,x_2)
\right)
u(a_1,\cosh(a_2-a_0)^{-1}x_1,z_1)\,v(a_2,\cosh(a_0-a_1)^{-1}x_2,z_2)
\times
$$

$$
\times
\exp\left(
\frac{-i}{\hbar}\left\{
\Omega(\sinh(2(a_2-a_0)),z_1)+\Omega(\sinh(2(a_0-a_1)),z_2)
+\Omega(\sinh(2(a_1-a_2)),z_0)
\right\}
\right)
\cosh(2(a_1-a_2))
$$
which, after changing the variables following 
$x_1\leftarrow\cosh(a_2-a_0)^{-1}x_1,\, 
x_2\leftarrow\cosh(a_0-a_1)^{-1}x_2$,
yields the announced formula.
\EPf

\begin{rmk}
{\rm 
It is important to mention that, in Theorem \ref{THM}, the function spaces $\e_\hbar$ are not invariant 
subspaces of $C_\infty(R)$ under the action of $R$. In order to obtain invariant spaces, 
one can consider completions of the $\e_\hbar$'s with respect to suitable $C^\star$-norms. This point has been 
treated in details in \cite{Bi2}.
}
\end{rmk}
\section{The case $n=1$}\label{Mozart}
In \cite{Bi2}, one finds a strict quantization of the symplectic symmetric space $M=(SO(1,1)\times\R^2)/\R$.
Let us first briefly recall how this quantization is defined. 
It turns out that the above mentioned symplectic symmetric space is, as a symplectic manifold, 
globally symplectomorphic to the two-dimensional symplectic vector space $(\R^2=\{(a,l)\}, da\wedge dl)$.
With respect to the coordinate system $(a,l)$, the geodesic symmetries 
are 
$$
s_{(a,l)}(a', l')=(2a-a',2\cosh(a-a')l-l').
$$
The curvature endomorphism of the underlying connection, $\nabla$, is given by 
$$
R(\partial_a,\partial_l)=\left(
\begin{array}{cc}
0&0\\
-1&0
\end{array}
\right).
$$
Defining the following family of diffeomorphisms~:
$$
\phi_\hbar:\R^2\to\R^2:\phi_\hbar(a,\alpha)=(a,\frac{2}{\hbar}\sinh(\frac{\hbar}{2}\alpha)),
$$
a result analogous to Lemma~\ref{PHI} allows us to define a linear injection
$$
\tau_\hbar:\s(\R^2)\to\s'(\R^2)
$$
by
$$
\tau_\hbar=F^{-1}\circ(\phi_\hbar^{-1})^\star\circ F,
$$
where
$$
Fu(a,\alpha)=\int_{\R} u(a,l)e^{-i\alpha l} dl.
$$
Entirely similarly as in Section~\ref{Bach}, the image space
$$
\e_\hbar=\tau_\hbar(\s(\R^2))\subset\s'(\R^2)
$$
is shown to be endowed with an associative product defined by
$$
a\star_\hbar b=\tau_\hbar(T_\hbar a\star^W_\hbar T_\hbar b),
$$
where $T_\hbar:\s(\R^2)\to\s(\R^2)$ is given by 
$T_\hbar=F^{-1}\circ\phi_\hbar^\star\circ F$, and 
where $\star^W_\hbar$ denotes the Weyl product on $\R^2$ (cf. Formula (\ref{WEYL})). A computation similar 
(but simpler) to the one in the proof of Theorem \ref{THM} leads us to the following formula
\begin{equation}\label{N=1}
u\star_\hbar v(x_0)=\frac{1}{\hbar^{4}}\int_{M\times M} 
e^{\frac{2i}{\hbar}S(x_0,x_1,x_2)}\cosh(a_1-a_2)\,u(x_1)\,
v(x_2) dx_1dx_2
\end{equation}
for $u,v\in\s(\R^2)\subset\e_\hbar$, where $x_i=(a_i,l_i)\in M=\R^2$, 
where $dx_i$ denotes the Liouville measure on $M$and where
$$
S(x_0,x_1,x_2)=\oint_{0,1,2}\sinh(a_0-a_1)l_2\quad(x_i=(a_i,l_i)\in M=\R^2).
$$
Besides associativity, the main property of the product $\star_\hbar$ 
is its invariance under the transvection group $G=SO(1,1)\times \R^2$ of the symmetric space $(M,\nabla)$. In other 
words, both amplitude $\cosh(a_1-a_2)$ and phase $S$ are invariant functions under the diagonal action of $G$.\\
Now, we observe that the transvection  group $SO(1,1)\times\R^2$ actually contains a subgroup $R$ isomorphic to the 
Iwasawa subgroup $AN$ of $SU(1,1)$. Indeed, the table of the Lie algebra $\g$ of $G$ is 
$$
\begin{array}{ccc}
\left[a,l\right] & = & k\\
\left[k,a\right] & = & -l\\
\left[k,l\right] & = & 0.
\end{array}
$$
Therefore, $\r=\mbox{span}\{a,k+l\}$ is a subalgebra isomorphic to $\a\times\n$ in $\mbox{su}(1,1)$. Observe that 
the analytic subgroup $R$ of $G$ with algebra $\r$ acts simply transitively 
on $M$. Hence, the quantization $\star_\hbar$ 
(cf. Formula (\ref{N=1})) defines a left-invariant strict quantization of the (symplectic) Lie group $R$. In 
particular, one can interpret Formula (\ref{N=1}) in two ways. One way is to say that it is a degeneracy of 
Formula (\ref{PRODUCT}) for a one-dimensional nilpotent factor $\n$. This 
emphasizes more the group representation theoretical 
aspect of the construction. The other way relies on the fact that the phase $S$ (as well as the amplitude) is 
determined uniquely in terms of the symmetric symplectic geometry of the symplectic symmetric space 
$(M,\omega,\nabla)$ \cite{W, Bi2}.\\
We end this section by mentioning an equivalence between our product formula (\ref{N=1}) in the degenerated case $n=1$ and 
Unterberger's formula for the composition of symbols in the one-dimensional 
Klein-Gordon Calculus (Formula (2.9) in \cite{U2}, 
see also \cite{U1}). More precisely, let $f_1$ and $f_2\in C_c^\infty(\R^2)$ be two compactly supported functions and let 
$f_1\sharp^U_\hbar f_2$ denote the symbol of the composition $Op(f_1)\circ Op(f_2)$ in the one-dimensional Klein-Gordon Calculus 
(see \cite{U2} pp. 174). Define the diffeomorphism 
$$
\varphi:\R^2\to\R^2:\varphi(a,l)=(\frac{1}{\cosh(a)}l,\sinh(a)).
$$
Then, one has 
$$
\varphi^\star(f_1\sharp^U_\hbar f_2)=(\varphi^\star f_1)\star_\hbar(\varphi^\star f_2),
$$
where $\star_\hbar$ is the product defined in Formula (\ref{N=1}).

\section{Remark for further developments}\label{Schubert}
Formulae (\ref{PRODUCT}) and (\ref{N=1}) define left invariant associative multiplications 
on the spaces $\e_\hbar$'s. 
The latter spaces play an analogous role the Schwartz space $\s(\r)$ does in the case of 
Weyl's quantization.
Each algebra $(\e_\hbar,\star_\hbar)$ is  isomorphic (via the ``equivalence'' $T_\hbar$) 
to $(\s(\r),\star^W_\hbar)$ 
(see Theorem \ref{THM}). It would therefore not be surprising that the deformed products 
(\ref{PRODUCT}) and (\ref{N=1}) 
extend to the space of smooth bounded functions, as in the case of Weyl's 
quantization \cite{R1}. 
This should provide actual universal 
deformations for any action of the group $R$ on any $C^\star$-algebra.

\end{document}